# Exact Anchoring and a Dualization-Based Matheuristic for Bi-Level Dual-Defense Network Interdiction

**Wei-Chang Yeh**
Department of Industrial Engineering and Engineering Management, National Tsing Hua University, Hsinchu 300, Taiwan, R.O.C.
E-mail: yeh@ieee.org

**Abstract –** Bi-level interdiction models are frequently solved by metaheuristics whose solution quality cannot be assessed, because exact optima are unavailable at the scales tested. We supply them for the bi-level dual-defense attacker model (BDAM), which couples node interdiction, edge destruction and capacitated supply support, and which we previously solved by a hybrid metaheuristic. First, BDAM's dominant attacker-path term admits an exact single-level reformulation by lower-level dualization, a reduction available whenever arc lengths are linear in the defender's binary decisions; the resulting mixed-integer program certifies optimality on all eighteen three-row and five-row configurations of our earlier benchmark, seventeen in under ten seconds, with strong uncertified incumbents out to 15×30 grids. Our published averages sit 3.55% below that frontier, and the gap widens with scale. Second, the supply rule, like any rule priced on a single attacker shortest path, is ill-posed under ties, and our threat-corridor formulation is tie-invariant by construction. The tie-break moves the objective by under $10^{-3}$ but the realized supply cost by up to 0.83 units, so the defect is suppressed by the objective weight rather than absent. Third, MILP-DA pairs the exact anchor with the corridor decode and a feasibility repair; compared against our published figures with no re-implementation on either side, it wins on thirteen of eighteen certified configurations and on all eighteen larger ones, every loss falling on a three-row grid.

## 1. INTRODUCTION AND MOTIVATION

Network interdiction studies how an upper-level interdictor should modify a network so as to worst-case-degrade the response of a lower-level user (Israeli & Wood, 2002; Smith & Song, 2020; Wood, 1993). The shortest-path network interdiction problem (SPNIP) is the canonical case: the interdictor maximizes the length of the shortest source-to-terminal ($s$-$t$) path that the user will then travel (Fulkerson & Harding, 1977; Israeli & Wood, 2002). A recent and expressive extension is the bi-level dual-defense attacker model (BDAM) of Yeh et al. (2025), which augments SPNIP with two coordinated defense operations, node interdiction and edge destruction, and with an explicit defender supply-support requirement, so that every interdicted node must be serviced from a capacitated supply node. BDAM is NP-hard and is solved by Yeh et al. (2025) with a hybrid metaheuristic, the improved simplified swarm optimization with Dijkstra (iSSOD). The present paper is a self-assessment of that model by one of its authors.

Expressiveness, however, comes at the cost of verifiability, and three specific weaknesses motivate this paper. First, BDAM provides no exact optimality anchor. Its reported solutions are compared only against other metaheuristics, so it is impossible to know how far a reported objective is from the true optimum, or whether the claim by Yeh et al. (2025) that conventional exact methods "fall short on large-scale problems" actually holds at the scales tested.

Second, the defender's supply routing is priced by the adjusted edge length of the single-path rule, Eq. (4) of Yeh et al. (2025), which we state as (3.4) in Section 3.2. There $y_{ij}$ records edge destruction, $p_{ij}$ records whether the attacker's chosen path traverses edge $(i, j)$, and $M$ is a large constant penalizing any supply route that overlaps that path. The difficulty is that p is not a defender decision: it solves the lower-level problem, and when several attacker paths tie at the optimum it is not unique. The supply cost $D_U$, and with it the objective, therefore depends on which optimal path the shortest-path routine happens to return rather than on the model itself, so the rule is strictly ill-posed. The single-path rule is the only point at which the attacker's solution enters the defender's cost, so any non-uniqueness in the lower-level optimum propagates directly into the upper-level objective. This weakness is not peculiar to BDAM: any lower-level rule priced on “the” shortest path inherits it.

Third, BDAM offers no guarantee that its supply requirement can be met. The randomized repair mechanism of Yeh et al. (2025) restores violations of the supply-capacity constraint and the destruction budget, but nothing ensures that an attacker-avoiding route to an assigned ward exists at all. Under the penalty pricing of (3.4) the difficulty is hidden, since a blocked route is charged M rather than rejected; under the implementation of Yeh et al. (2025), which removes the attacker's edges outright, a ward can simply be unreachable and the objective collapses to the penalty. The defender's logistical requirement, which is the model's stated advance over conventional SPNIP, is thus guaranteed in neither reading. None of the three weaknesses is fatal to BDAM as an operational model, but all three must be resolved before its solutions can be trusted as optimal or its algorithmic claims accepted.

This paper resolves all three, and makes three corresponding contributions.

1. **Exact anchors.** We show that the dominant (attacker-path) term of BDAM admits an exact single-level reformulation by lower-level dualization in the style of Israeli and Wood (2002). Because the interdicted arc length is linear in the defender variables, maximizing the attacker's shortest path becomes a compact MILP requiring no big-M linearization, solvable by the free CBC solver; the same reduction is available for any interdiction model whose arc lengths are linear in the defender's binary decisions. The eighteen three-row and five-row configurations of our earlier benchmark are solved to certified optimality, giving, apparently for the first time for this model family, a set of exact reference optima $H^*$. For the eighteen larger configurations the same formulation returns strong incumbents with certified dual bounds, but no proof of optimality within the time budget.
2. **A well-posed supply rule.** We make the single-path rule well-posed by blocking, for supply routing, not one arbitrary shortest path but every edge lying on some attacker shortest path: the edges of the attacker's shortest-path directed acyclic graph, which we call the threat corridor. This set is determined by the instance alone and is therefore invariant to tie-breaking, and it is the operationally natural reading of the rule, since supply should avoid the whole region the attacker may traverse rather than one representative route. The single-path rule is genuinely ill-posed in integer-weight regimes: a bounded constructed instance exhibits a 28.6% swing in realized supply cost, and a Manhattan instance flips feasibility entirely, purely from the tie-break. On our own $\sqrt{2}$-diagonal instance family the effect is dormant in the combined objective, whose spread over twenty tie-breakings stays below $8.3\times10^{-4}$, but live in the supply cost beneath it, which the tie-break moves by up to 0.83 units. What suppresses it is the objective's weighting of the supply term, not the supply routing itself.

3. **MILP-DA matheuristic.** We pair the exact dominant-term solution of the first contribution with the threat-corridor decode of the second and a supply-feasibility repair, yielding MILP-DA: a full BDAM solution carrying an exactness certificate on the dominant term. Compared directly against the results published by Yeh et al. (2025), with no re-implementation of the metaheuristic on either side, MILP-DA's objective exceeds the published best of thirty runs on thirteen of the eighteen certified configurations (one-sided Wilcoxon signed-rank, $p=0.0035$) and on all eighteen larger ones ($p=1.1\times10^{-4}$), and exceeds the published averages on sixteen of eighteen. Every loss occurs on a three-row grid and none exceeds 0.0060 in objective value; MILP-DA wins all 27 configurations with five or more rows. The repair closes the third weakness: every ward in the returned plan is supplied under the corridor rule, by construction rather than at a rate determined by the tie-break.

Section 6 develops the resulting position, which cuts both ways. Certification proves cheap exactly in the regime where the metaheuristic was benchmarked, where our own published averages already leave a measurable gap, and breaks beyond it, so the large-scale justification for a metaheuristic retains genuine merit while its optimality at small and medium scale does not survive contact with an exact anchor. A second finding cuts across the first: the supply term that the dominant-term anchor does not capture grows with scale, so on the largest grids the restriction to the dominant term costs more than the missing certificate does.

The remainder of the paper is organized as follows. Section 2 reviews related work. Section 3 recalls BDAM, states the single-path rule (3.4) and analyzes its well-posedness, develops the shortest-path-DAG theory used throughout, and gives the tie-invariant threat-corridor formulation. Section 4 presents the exact anchor and the MILP-DA matheuristic. Section 5 reports the computational study on our earlier benchmark. Section 6 states the limitations, and Section 7 concludes.

## 2. RELATED WORK

This section positions the paper along two lines: exact interdiction methods, and the BDAM model and its metaheuristic solver.

### 2.1 Exact interdiction methods

Deterministic network interdiction originates with Wood (1993). Israeli and Wood (2002) formalized shortest-path interdiction as a bi-level max-min problem in which the interdictor adds a delay $d_a$ to the nominal length $c_a$ of each interdicted arc a under a resource budget, and the attacker then travels the cheapest *s-t* path:

$$max_{x \in X}\, min_{y \in Y}\ \sum_{a \in A}( c_a + d_a x_a) y_a \tag{2.1}$$

where $x_a \in \{0,1\}$ is the interdiction decision, $r_a > 0$ is the resource consumed by interdicting arc $a$, $B$ is the interdictor's total resource budget, $X = \{x: \sum_a r_a\, x_a \leq B\}$, and $y$ encodes an *s-t* path through the flow-conservation polytope. The inner shortest-path linear program can be dualized with node potentials $\pi\,(\pi_s$ =0), turning the bi-level problem into a single-level mixed-integer program:

$$max_{x,\pi}\, \pi_t\ \text{ s.t. } \pi_j - \pi_i \leq c_{ij} + d_{ij} x_{ij}\ \ \forall (i,j) \in A,\ \sum_a r_a\, x_a \leq B \tag{2.2}$$

Because the interdicted length is linear in $x$, (2.2) is a plain MILP requiring no big-M linearization, and its dual potentials furnish a per-instance optimality certificate; this is exactly the structure our anchor exploits in Section 4.1.

Israeli and Wood (2002) also give an equivalent covering/decomposition reformulation solved by attacker-path (row) generation, and Smith and Song (2020) survey the resulting exact and decomposition algorithms comprehensively. In the continuous-budget setting, Fulkerson and Harding (1977) maximize the minimum s-t path length by a parametric max-flow argument whose optimum equalizes a family of tied shortest paths across a minimum cut, a min-cut/shortest-path-DAG duality that reappears in our shortest-path-DAG theory (Section 3.4).

Most-vital-arc and most-vital-node analysis (Ball et al., 1989; Corley & Sha, 1982) identifies the components whose removal most degrades *s-t* performance, $v^* \in \arg\max_{v} d_{G-v}(s,t)$; the budgeted set version is NP-hard (Ball et al., 1989; Smith & Song, 2020). Resilience-oriented interdiction and critical-component identification have likewise been studied in infrastructure settings (Alderson et al., 2018; Bruneau et al., 2003; Church et al., 2004; Sharkey et al., 2021). These exact and combinatorial results anchor the dualization approach of this paper.

Although each component of BDAM is individually well studied, we are not aware of prior work combining node interdiction, edge destruction and a capacitated defender supply requirement in a single bi-level model, nor of any exact treatment of such a combination.

### 2.2 BDAM and metaheuristic interdiction

BDAM (Yeh et al., 2025) extends shortest-path interdiction with dual defense operations (node interdiction and edge destruction) and with defender supply support, and is solved by iSSOD, a population-based metaheuristic in which the attacker's response is evaluated by Dijkstra's algorithm (Dijkstra, 1959) inside the fitness function. Its core operator is the simplified-swarm stepwise update, which for each decision variable copies the global best $g_j$, the personal best $p_{ij}$, the current value, or a random value according to fixed cumulative thresholds:

$$x_{ij}^{(t)} = \begin{cases} g_j, & \rho < C_g \\ p_{ij}, & C_g \leq \rho < C_p \\ x_{ij}^{(t-1)}, & C_p \leq \rho < C_w \\ \text{rand}, & \rho \geq C_w \end{cases} \tag{2.3}$$

with $(C_g, C_p, C_w)$=(0.65,0.85,0.99) in Yeh et al. (2025), so a variable copies the global best with probability 0.65 and randomizes with probability only 0.01; this near-absent random slice is why Yeh et al. (2025) add an entropy-guided dynamic-swap local search for exploration. Metaheuristics are attractive for large, tightly-coupled bi-level models but provide no optimality certificate; our contribution is to supply the missing exact reference frame rather than to replace the metaheuristic.

## 3. THE BDAM MODEL AND THE WELL-POSEDNESS OF THE SINGLE-PATH RULE

This section sets out the model and the two properties of it that the rest of the paper depends on. Section 3.1 recalls BDAM. Section 3.2 states the single-path supply rule (3.4) in full, since the second contribution repairs it, and settles the convention under which a blocked supply route is read. Section 3.3 shows that the rule is ill-posed under shortest-path ties. Section 3.4 develops the shortest-path-DAG theory used throughout, and Section 3.5 uses it to give a tie-invariant threat-corridor formulation and measures what the tie-break is worth on our own instance family, where the effect proves dormant in the objective and live in the supply cost beneath it.

**3.1 Model overview**

BDAM is defined on an undirected network *G*=(*V*, *E*) with positive edge lengths $d(u,v)$, a source *s*, a terminal *t*, an eligible interdiction-node set, an edge-destruction budget *R*, and a set of capacitated supply nodes. The defender chooses (i) which eligible nodes to interdict, encoded by binary $X_v$, (ii) which eligible edges to destroy, encoded by binary $y(u,v)$ with total at most *R*, and (iii) a capacitated assignment of each interdicted node (ward) to a supply node. Interdiction and destruction lengthen arcs for the attacker: the interdicted arc length is

$$\delta^A(u,v) = d(u,v)(1 + g_1 X_v + g_2 y(u,v)) \quad (3.1)$$

where $g_1$ and $g_2$ weight node interdiction and edge destruction. The attacker then travels the shortest *s*-*t* path of length

$$H = min_{P \in \mathcal{P}_{st}} \sum_{(u,v) \in P} \delta^A (u,v) \quad (3.2)$$

The defender additionally routes supply to each ward while avoiding the attacker's path, incurring a supply cost $D_U$, and maximizes the combined objective

$$F = H - wD_U, w = 0.001 \quad (3.3)$$

**3.2 The single-path supply rule and the exclusion convention**

The rule by which that supply avoids the attacker's path is the object of this paper's well-posedness analysis, so we state it here rather than by reference. It is Eq. (4) of Yeh et al. (2025), reproduced here in this paper's notation (subscripts written as arguments, and their $\gamma_2$ written $g_2$ as in (3.1)), and we call it the single-path rule: the defender's adjusted length of edge $i,j$ is

$$\delta^D(i,j) \; = \; d(i,j)\big(1 + g_2\, y(i,j)\big) \; + \; M\big(p(i,j) + p(j,i)\big), \forall\, i,j \in V, \quad (3.4)$$

where $p(i,j)$ indicates that the attacker's selected path traverses edge $(i,j)$ and $p(j,i)$ its reverse, $y(i,j)$ is the edge-destruction decision, and $M$ is a large constant penalising any supply route that reuses an edge of the attacker's path. The defender's total supply cost aggregates this over the |*S*| supply routes (Yeh et al., 2025, Eq. 5):

$$D_U \; = \; \sum_{k=1}^{|S|} \; \sum_{(i,j) \in E} \delta^D\,(i,j)\, z(i,j,k), \quad (3.5)$$

with $z(i,j,k) = 1$ when supply route k uses edge $i,j$.

Note what the $M$ term does in (3.4): it prices a supply route that reuses an attacker edge rather than forbidding it, so $D_U$ as defined by (3.5) is always finite and no supply plan is ever infeasible under the original formulation of Yeh et al. (2025), the route simply costing at least $M$. The implementation described by Yeh et al. (2025, Section 4.3) does something different: it assigns infinite length to the attacker's edges, removing them from the defender's supply planning altogether. Model and implementation therefore diverge, and the ill-posedness established in Section 3.3 holds under both readings, appearing as an $M$-sized collapse of the objective under the first and as outright infeasibility under the second.

The threat-corridor formulation of Section 3.5 departs from (3.4) deliberately and follows the second reading: it excludes the blocked edges instead of pricing them, and a ward for which no corridor-avoiding route then exists is genuinely unsupplied. Our implementation charges $M = 10^6$ per such ward, a value Yeh et al. (2025) leave unspecified. That is the origin of the $2.0\times10^6$ entry in Table 2: two unsupplied wards, which at $w = 0.001$ drives $F$ to −1992.38. Every count of "supply-infeasible" runs in this paper is therefore a count under exclusion; the same situation under (3.4) appears as a finite route cost of at least $M$.

Because $w$ is small, $H$ is the dominant term; $D_U$ acts as a tie-breaker among interdiction plans of equal attacker delay. This separation of scales is central: it is exactly what makes an exact anchor on $H$ both possible and informative (Section 4.1).

**3.3 The single-path rule is ill-posed under shortest-path ties**

The rule (3.4) charges a large penalty $M$ for any supply route that uses an edge of the attacker's shortest path $P$. When several shortest paths tie, $P$ is selected by the tie-breaking order of the shortest-path routine, so the blocked edge set, and hence the supply cost $D_U$ and the objective $F$, depends on an implementation detail rather than on the model. We make this precise with a bounded, always-feasible integer-weight instance (Fig. 1): two attacker paths tie at $H$=10, and the defender's ward can always be supplied, but the cheapest supply route reuses an edge of the southern attacker path. If the tie-break selects the southern path, that edge is blocked and the supply detour costs more.

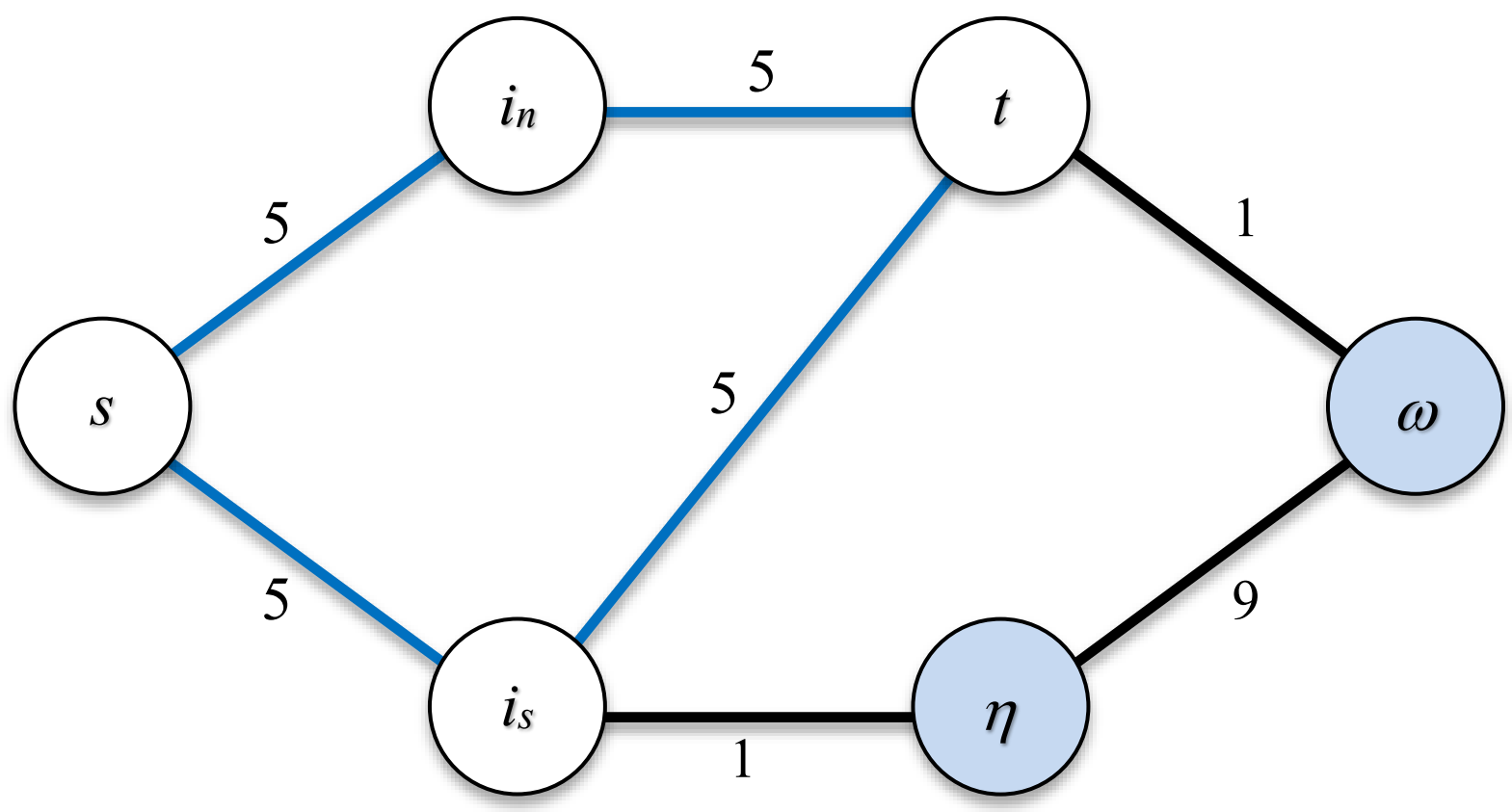

**Fig. 1.** A bounded, always-feasible integer-weight instance. The two attacker $s$-$t$ routes (blue) tie at $H$ = 10; the ward $\omega$ is supplied from hub $\eta$ (grey). The cheapest supply route reuses edge $i_s$-$t$, so the tie-break that selects the southern attacker path raises the supply cost from 7 to 9 (a 28.6% swing) with $H$ unchanged.

The instance fixes the interdiction plan and exhibits the supply rule in isolation, so $\omega$ is a ward already designated by some $X$ and the tie-break is the only quantity that varies.

Sampling sixty seeded draws of the node-priority order used as the shortest-path routine's tie-breaker (seeds 0–59, which yield 56 distinct orderings), the realized supply cost swings between 7 and 9 units, taking the value 7 on thirty of the sixty draws and 9 on the other thirty, a 28.6% change driven entirely by the arbitrary tie-break, with $H$ unchanged. On a pure Manhattan (integer, no-diagonal) grid, where ties are pervasive, the effect is worse: some tie-breaks block the ward's only corridor, flipping the solution from feasible to infeasible and collapsing the objective. Thus the single-path rule is genuinely ill-posed in integer-weight regimes.

### 3.4 Shortest-path-DAG theory

Both the corridor formulation and the matheuristic rely on the structure of the set of shortest paths. Let $d_s(v)$ and $d_t(v)$ be the shortest distances from $s$ and to $t$ under a fixed length function, and let $L_0 = d_s(t)$. An undirected edge $(u, v)$ lies on some shortest $s$-$t$ path iff

$$d_s(u) + d(u, v) + d_t(v) = L_0 \text{ or } d_s(v) + d(u, v) + d_t(u) = L_0 \qquad (3.6)$$

Orienting each such edge in the direction that satisfies the equality yields the shortest-path DAG, $G_0$. The following corrects a subtlety in path-based enumeration: the shortest-path subgraph must be oriented, because an unoriented subgraph can admit walks that are not shortest paths.

**Theorem 1 (oriented path-closure).** The directed *s-t* paths of the oriented SP-DAG $G_0$ are exactly the shortest *s-t* paths of *G*.

**Proof.** (Soundness) If $(u, v)$ is oriented $u \to v$ then $d_s(u) + d(u, v) + d_t(v) = L_0$, so $d_s(v) = d_s(u) + d(u, v)$; summing telescopically along any directed *s-t* path gives total length $d_s(t) = L_0$, a shortest path. Because $d_s$ strictly increases along every arc, no directed path repeats a vertex, so it is simple. (Completeness) If *P* is a shortest path, then for each traversed edge the *s*-to-*u* prefix, the edge, and the *v*-to-*t* suffix sum to $L_0$, so the edge satisfies the equality in its traversed direction and is present in $G_0$ with that orientation; hence *P* is a directed s-t path of $G_0$. ■

The orientation is what the earlier undirected statement omitted: without it, an edge satisfying the equality in one direction could be traversed in the other, producing a non-shortest or non-simple walk. Orientation also gives the mixing property.

**Proposition 1 (path mixing, $p_{mix}$).** If two shortest paths share an intermediate node *m*, then the path formed by following the first up to *m* and the second after *m* is also a shortest path. Consequently, if the SP-DAG has a node reachable by *a* prefixes and leading to *b* suffixes, at least *a·b* distinct shortest paths exist, so $|P^*|$ can be exponential in $|V|$.

**Proof.** Concatenation at *m* preserves the arc equalities on both segments, so the mixed walk has length $L_0$ and, by monotonicity of $d_s$, is simple; it is therefore a shortest path. The counting claim follows by choosing prefixes and suffixes independently. ■

This is why enumerating all shortest paths is output-sensitive and why interdiction based on a single shortest path is unsafe: mixing can restore an *s-t* route that a single-path analysis believed destroyed.

The six-node network of Fig. 2 serves as the running example for the orientation rule and illustrates mechanisms that the four-node example of Yeh et al. (2025) cannot. Nodes are $\{s, a, b, c, d, t\}$; unit edges form three routes: *s-a-t*, *s-b-t*, and *s-c-d-t*. The two two-hop routes tie at $L_0 = 2$ (a genuine shortest-path tie), while *s-c-d-t* has length 3.

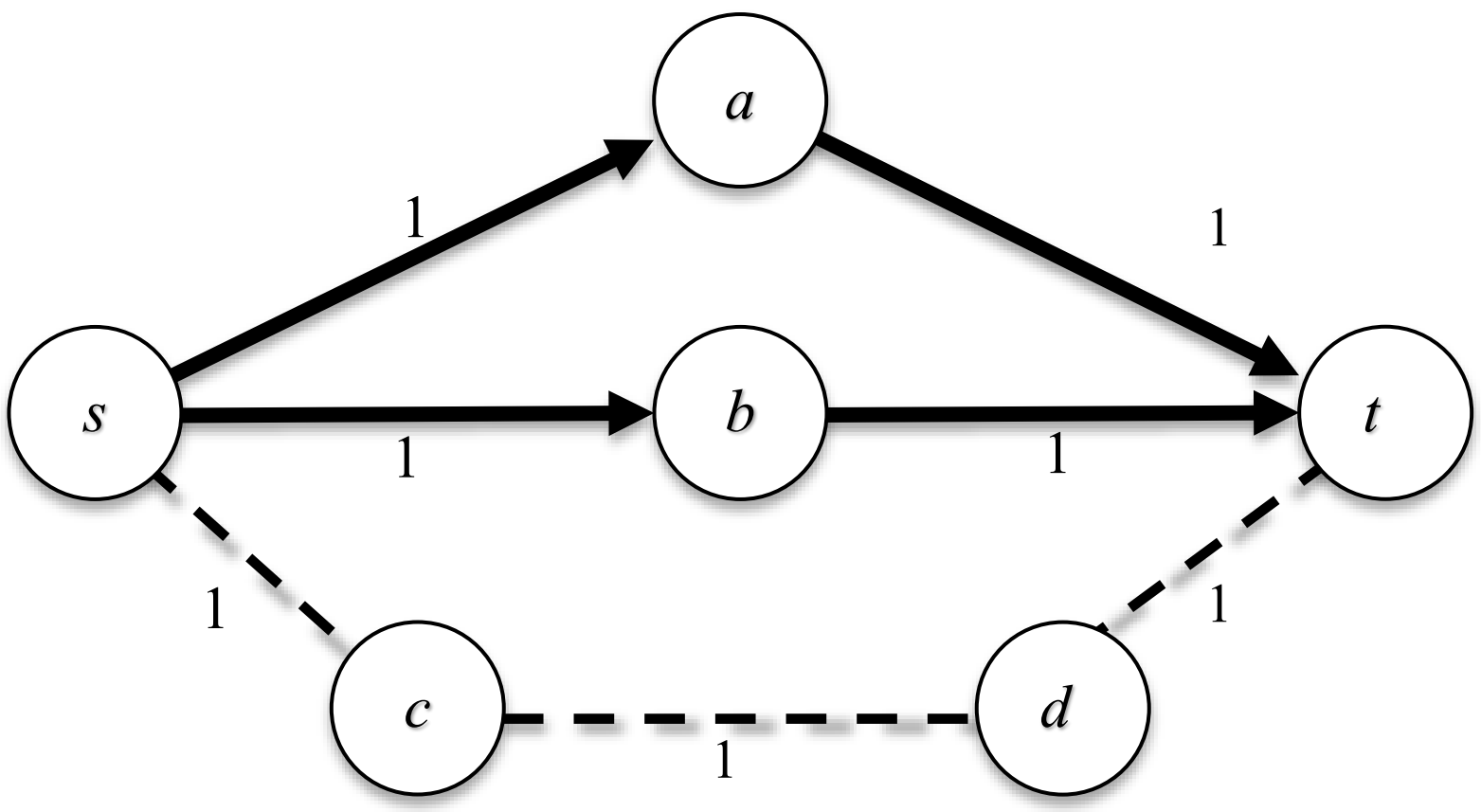

**Fig. 2.** The running six-node example (unit edges). The shortest paths *s–a–t* and *s–b–t* tie at $L_0 = 2$ and share no intermediate node, so the common intermediate set is empty and no single-node interdiction can lengthen the attacker's route. Heavy arrows mark the edges retained and oriented on the shortest-path DAG; the dashed edges *s–c*, *c–d* and *d–t* are off-DAG detours.

Enumerating shortest paths returns exactly $\{s\text{-}a\text{-}t,\ s\text{-}b\text{-}t\}$; their intermediate-node sets are $\{a\}$ and $\{b\}$, so the common intersection is empty.

Table 1 orients the edges of the six-node graph by the rule just stated: an edge $(u, v)$ lies on the SP-DAG iff $d_s(u) + d(u, v) + d_t(v) = L_0$ with $L_0 = 2$, and is then oriented toward increasing $d_s$. Four edges are retained and oriented $s \to a \to t$ and $s \to b \to t$; the three edges *s-c*, *c-d*, *d-t* are off-DAG detours. The two oriented *s-t* paths are exactly the shortest paths, and because they share no intermediate node the common intersection is empty, so no single-node interdiction can degrade service. The supply-feasibility repair cannot be shown here, because this example carries no supply term; it is traced separately in Section 4.2.

**Table 1.** SP-DAG orientation for the six-node example[a]

| Edge $(u, v)$ | $d_s(u)$ | $d_t(v)$ | On SP-DAG? | Orientation |
|---|---|---|---|---|
| *s-a* | 0 | 1 | yes | $s \to a$ |
| *s-b* | 0 | 1 | yes | $s \to b$ |
| *a-t* | 1 | 0 | yes | $a \to t$ |
| *b-t* | 1 | 0 | yes | $b \to t$ |
| *s-c* | 0 | 2 | no | off-DAG (detour) |
| *c-d* | 1 | 1 | no | off-DAG (detour) |
| *d-t* | 2 | 0 | no | off-DAG (detour) |

[a] $L_0 = 2$; edges satisfying (3.6) are retained and oriented.

Proposition 1 is equally concrete on this graph: the layered variant $s \to \{a, b\} \to m \to \{c, d\} \to t$ has $2{\times}2 = 4$ shortest paths, every one of which passes through the single mixing node m, so interdicting $\{m\}$ lengthens all four at once, whereas interdicting any peripheral node leaves three of them at their original length.

### 3.5 A well-posed threat-corridor formulation

The remedy is to block, for supply routing, not one arbitrary shortest path but every edge that lies on some attacker shortest path. We call that set the threat corridor and write it, for a defender plan $X, y$ inducing the arc lengths $\delta^A$ of (3.1),

$$C(X,y) \;=\; \{\,(u,v) \in E \,:\, (u,v) \text{ satisfies (3.6) under } \delta^A\,\}, \tag{3.7}$$

that is, the edge set of the shortest-path DAG $G_0$ of Section 3.4. The single-path rule (3.4) is then replaced by the corridor rule, which excludes rather than prices:

$$\delta_C^D(i,j) \;=\; \begin{cases} d(i,j)\big(1 + g_2\, y(i,j)\big), & (i,j) \notin C(X,y), \\ \infty, & (i,j) \in C(X,y), \end{cases} \tag{3.8}$$

with the supply cost aggregated as in (3.5) and a ward admitting no finite route charged the infeasibility penalty $M = 10^6$ of Section 3.2. Blocking the corridor is also the operationally natural reading of the rule: supply should avoid the whole region the attacker may traverse, not a single representative route. The formulation is well posed, which is the property (3.4) lacks.

**Proposition 2 (tie-invariance of the corridor).** For a fixed instance and a fixed defender plan $X, y$, the set $C(X,y)$ of (3.7) is uniquely determined, and $\delta_C^D$ and the resulting supply cost $D_U$ are therefore independent of any tie-breaking among the attacker's optimal paths.

**Proof.** Under $\delta^A$ the distance functions $d_s$ and $d_t$ are the unique shortest-distance functions of the network, and $L_0 = d_s(t)$ is likewise unique; none of the three depends on which optimal path a shortest-path routine returns. Membership in $C(X,y)$ is decided by (3.6), a test referring only to $d_s$, $d_t$, $d$, and $L_0$. Hence $C(X,y)$ is a function of the instance and $X, y$ alone. By contrast, the blocked set $E(P)$ of (3.4) is defined from a selected optimal path $P$, which is not unique when ties occur. ■

**Corollary 1 (the corridor is conservative).** $E(P) \subseteq C(X,y)$ for every attacker-optimal path $P$. Consequently the corridor rule blocks a weak superset of the edges blocked by (3.4), its supply cost is weakly larger, and the resulting objective $F = H - wD_U$ is weakly smaller.

**Proof.** By the completeness half of Theorem 1, every edge of a shortest path satisfies (3.6) in its traversed direction and so lies in $C(X,y)$. The cost and objective statements follow from the monotonicity of (3.5) in the blocked set. ■

Corollary 1 has a practical consequence recorded here and quantified in Section 5.3: because MILP-DA is scored under the corridor while the figures published by Yeh et al. (2025) are scored under the single-path rule, every comparison in Section 5 handicaps MILP-DA, and its reported wins therefore hold a fortiori. Well-posedness matters in theory; we tested whether it matters on our own instances. Re-evaluating the best found solutions under twenty randomized tie-breakings on the √2-diagonal grid family with the weight w = 0.001, the single-path objective spread never exceeds $8.3\times10^{-4}$, because the diagonal lengths make exact ties rare and w makes the supply term second-order. Under the corridor rule the spread is exactly zero on every row, at no measurable cost. The effect is nevertheless dormant in F rather than absent: the same tie-breakings move the realized supply cost $D_U$ itself by up to 0.83 units (instance 4, 3×5: $D_U$ ranges 15.66–16.49, or 5.3% of its smallest value), and it moves on 3 of the 6 instances tested. What suppresses the swing is the weight, not the stability of the supply routing, so any use of the model that reads $D_U$ directly, through a supply budget, a service-level constraint, or any larger *w*, sees it in full. Section 5.7 reports the measurement in detail (Table 8).

## 4. ALGORITHMS

Section 3 leaves two things in place: a dominant term that is linear in the defender's decisions, and a supply rule that is well posed once the threat corridor replaces the single path. This section turns each into an algorithm. Section 4.1 exploits the linearity of (3.1) to obtain an exact anchor for the dominant term by lower-level dualization; Section 4.2 pairs that anchor with the corridor decode of Section 3.5 and a supply-feasibility repair, giving MILP-DA, a full BDAM solution that carries a dominant-term certificate.

### 4.1 Exact dominant-term anchor by dualization

The attacker's problem is a shortest-path linear program; its dual is a longest-potential program. Specializing the generic dualization (2.2) to BDAM's linear arc length (3.1) and introducing a node potential $\pi_v$, the dominant term equals

$$H = max\,\pi_t \ \text{s.t.}\ \pi_s = 0,\ \pi_v - \pi_u \le \delta^A(u,v)\ \forall(u,v) \tag{4.1}$$

Because $\delta^A(u,v) = d(u,v)\big(1 + g_1 X_v + g_2\, y(u,v)\big)$ is linear in the binary defender variables, maximizing $H$ over the defender's decisions, subject to the capacitated ward-to-supply assignment and the edge-destruction budget $\sum_{(u,v)\in E} y\,(u,v) \le R$, is a single-level MILP. We solve it with the open-source CBC solver (Forrest & Lougee-Heimer, 2005).

The optimum $H^*$ is an exact anchor for the dominant term, because $F = H - wD_U$ with $0 \leq wD_U \leq wD_U^{\max}$, where $D_U^{\max}$ is any finite upper bound on the supply cost, which exists because under the exclusion convention of Section 3.2 every ward is either served at finite cost or charged the penalty $M$. It follows that $H^* - wD_U^{\max} \leq F^* \leq H^*$, and $\text{gap(alg)} = \frac{H^* - H_{\text{alg}}}{H^*}$ is a valid optimality gap for the term that determines BDAM's objective to first order.

**4.2 MILP-DA: a dualization-based matheuristic**

MILP-DA turns the anchor into a full, deployable BDAM solution: it solves the dualization MILP (4.1) for an optimal interdiction and destruction plan, selects among the plans attaining the same $H^*$ the one of least supply cost, decodes the supply routes under the well-posed threat-corridor rule of Section 3.5, and repairs any ward left unsupplied. Algorithm 1 states this precisely. The dominant-term optimality certificate is inherited from the first solve and is unaffected by the selection, which only chooses among plans already attaining $H^*$.

The selection stage is not cosmetic. The equal-objective set has a dimension beyond the interdiction and destruction decisions: two plans may share both $X$ and $y$ yet assign wards to different supply nodes, which changes $D_U$ while leaving $H$ untouched. Without that stage, STEP 1 returns an arbitrary member of the set and the reported $F$ depends on which optimum the solver happened to reach first. STEP 2 searches that set, so plan selection is determined by the model rather than by the solver, and the retained plan is persisted so that the reported F is reproducible from the package rather than resampled at each solve.

The repair is likewise necessary. The MILP maximizes H and is blind to whether a ward can actually be supplied once the threat corridor is blocked; a geometrically valid interdiction plan can therefore leave a ward whose corridor-avoiding supply route does not exist. This is not hypothetical: on the reconstructed $r$=3, $c$=4 grid the raw MILP solution is supply-infeasible, as it is on several reconstructed $r$=5, $c$=5 instances, so a ward must be released to restore feasibility. The repair iteratively drops the offending ward, the least-value interdiction whose supply route is blocked, and re-evaluates until every remaining ward is supplied. Releasing a ward removes its interdiction and can lower $H$; MILP-DA thus trades a small, quantified amount of dominant-term optimality for guaranteed supply feasibility, a trade the raw metaheuristic makes silently and without a certificate.

**Algorithm 1:** MILP-DA

**INPUT:** BDAM instance; weights $g_1$, $g_2$; edge budget $R$; supply capacities.

**STEP 1.** Solve the dualization MILP (4.1) for an optimal $X$, $y$ and ward-to-supply assignment $I$; record the exact $H^*$ certificate.

**STEP 2.** Holding the attacker potential at $H^*$, search alternative optima of (4.1) and retain the plan of least supply cost $D_U$; persist the retained plan ($X$, $y$, $I$). The search is deterministic given its parameters but is not proven exhaustive, so the retained $D_U$ is an upper bound on the true minimum and the returned $F$ is a valid lower bound (Section 5.3).

**STEP 3.** Build the threat corridor (SP-DAG edges, Section 3.4) and route each ward's supply on the corridor-avoiding network; compute $D_U$.

**STEP 4.** While some ward is unsupplied ($D_U$ hits the infeasibility penalty): drop that ward from $I$ and re-evaluate $D_U$.

**STEP 5.** Return the full plan ($X$, $y$, $I$), the realized $F = H - w \cdot D_U$, and the $H^*$ certificate.

One point about STEP 4 belongs in the text rather than in the pseudocode. The repair exists because of the exclusion convention: under (3.4)'s pricing a ward is never unsupplied, so no repair would be needed and no certificate of feasibility would be meaningful; under exclusion it can be, and the repair is what restores feasibility at a quantified cost in $H$. Searching the assignment dimension alone, without varying the interdiction plan, improves 2 of the eighteen certified instances by 1.000 and 1.414 units of supply cost with $H$ unchanged on all eighteen (drift 0.00e+00), a real but small effect that STEP 2 largely already captures.

The repair is best seen on a genuine BDAM grid, since the six-node example of Section 3.4 carries no supply term. It must be shown on a reconstructed instance: on the authors' own benchmark the repair never fires, MILP-DA's decoded plan being supply-feasible on all 36 instances (dropped = 0 throughout, Section 5.3), and the 3×4 geometry used below is one the real family does not contain (Section 5.1).

On the reconstructed ($r$=3, $c$=4) instance (seed 304, certified $H^*$ = 7.6213) the dominant-term MILP interdicts four wards, but under the well-posed threat corridor two of them cannot be supplied, so the raw objective collapses to the two-ward infeasibility penalty of $2.0\times10^6$ (Table 2, decoded). The repair drops the least-value blocked ward, grid node (0, 2); releasing it removes its interdiction, which re-shapes the threat corridor, and on re-evaluation the other blocked ward, grid node (1, 3), is supplied, so a single pass restores feasibility at a small, quantified cost in the dominant term ($H$ falls 7.6213 → 7.1213), and returns $F$ = 7.1130. Two counterfactual decodes isolate the mechanism. Releasing node (0, 2)'s interdiction so that the corridor is recomputed, with every other ward's assignment left as it stood at iteration 0, restores ward (1, 3) at a supply distance of 3.50; holding the corridor at its iteration-0 state while releasing node (0, 2) from supply routing leaves ward (1, 3) blocked. The corridor re-shape is therefore decisive and the freed assignment is not a second mechanism. The decode prices each ward's route as an independent shortest path from its assigned supply node, with no capacity coupling between wards, so a released assignment cannot by itself restore another ward's reachability; that is what makes the repair a structural effect rather than an accounting one. This is exactly the H-for-feasibility trade the metaheuristic makes silently and without a certificate.

**Table 2. Supply-feasibility repair trace on the reconstructed (3,4) instance, ($g_1$, $g_2$) = (0.5, 0.5)[a]**

| Repair iteration | Action | wards \|*I*\| | *H* | $D_U$ | *F* | feasible? |
|---|---|---|---|---|---|---|
| 0 | raw MILP solution | 4 | 7.6213 | $2.0\times10^6$ | −1992.38 | no (penalty) |
| 1 | drop ward (0,2) | 3 | 7.1213 | 8.33 | 7.1130 | yes |

[a] seed 304, threat-corridor referee. Iteration 0 is the plan as returned by STEP 3, before any repair; iteration 1 is after one pass of STEP 4. [family=reconstructed; source_json=res_milpda_full.json; time_limit=1200 s; concurrency=solo] [run_id=n/a; pre-dates run stamping]

## 5. COMPUTATIONAL STUDY

This section reports the study on our own thirty-six benchmark instances. Section 5.1 describes that family and the documented reconstruction used alongside it. Sections 5.2 to 5.4 carry the paper's claims and rest entirely on figures printed by Yeh et al. (2025), so no re-implementation of the metaheuristic enters them: Section 5.2 measures our published averages against the certified anchors, Section 5.3 compares MILP-DA against the published columns instance by instance, and Section 5.4 states the tests and the family of record. Section 5.5 reports the reconstruction in full and explains why nothing rests on it. Sections 5.6 and 5.7 locate the certified-optimality frontier and check tie-sensitivity.

**5.1 Instances: the authors' benchmark and a documented reconstruction**

For the exact-anchor comparison of Section 5.2 we use the authors' own benchmark instance files (afprob2025.mat), which contain all 36 problems of Yeh et al. (2025). We verified their structure: an $r \times c$ grid with √2 diagonals, the attacker source and terminal, the right-edge supply nodes, the protected (non-interdictable) zone around $s$ and $t$, and a per-supply capacity $c_j = \lceil 0.2|V|/|S| \rceil$ (stored in the field nItd; validated on instances 1, 10, 19, 34 as 2, 3, 5, 13). The node-interdiction capacity is therefore $|S| \cdot c_j$ and the edge-destruction budget is $R$ = nDtr. Independently, we also regenerated the same grid family with a documented generator (assumptions [A1]–[A5]); it reproduces the real instances' grid backbone, source, terminal, and eligible-node rule exactly, differing only in the random diagonal draw [A1].

Sections 5.2–5.4 report results on the real instances throughout; the reconstruction is retained as an independent cross-check of the generator and of the supply decode (on all eighteen reconstructed configurations the shared decode reproduces the original one's plan and objective to nine decimals), and it is the family on which the earlier reconstructed-instance statistics quoted in Section 5.5 were computed. One difference between the two families proved consequential: a language-level rounding convention (Python's round-half-even: round(4.5)=4) silently substituted a 3×4 geometry for the benchmark's 3×5, and that artifact geometry alone accounted for 67% of the reconstructed supply-infeasibility events (Section 5.5).

Instances 19–36 are reported from the original stretch run, whose plans were not persisted, and their anchors are best-effort incumbents rather than certified optima. A fresh solve does not necessarily re-attain them: re-solving instance 19 under the same protocol returns $H$=19.364 against an anchor of record of 20.107, a shortfall of 3.7%, larger than the mean quality-interval width of 2.96% on that tier. The anchors remain valid lower bounds on the true optimum, so every comparison against them is conservative and no claim weakens.

### 5.2 Certified anchors versus the published metaheuristic results

On the real instances we solve the dominant-term anchor $H^*$ by lower-level dualization under the capacity model (per-supply $c_j$, edge budget $R$) with the CBC solver, and compare it against the metaheuristic's own published average objective $F_{\text{avg}}$ from Yeh et al. (2025, Tables 9–10), requiring no re-run of the metaheuristic.

Because $F_{\text{avg}} \le F^* \le H^*$ (an achieved average cannot exceed the true optimum, which cannot exceed the dominant-term optimum), the quantity $(H^* - F_{\text{avg}})/H^*$ bounds how far the published results sit below the certified frontier. Table 3 covers the $r$=3 and $r$=5 families (instances 1–18; the $r$=10 and $r$=15 families are anchored in Section 5.6). The published averages leave a mean gap of 3.55% (min 0.13%, max 7.52%), and the gap grows sharply with scale: 1.40% at $r$=3 versus 5.70% at $r$=5.

Thus, measured against certified optima on the authors' own instances and using the authors' own reported averages, the published metaheuristic is demonstrably sub-optimal, and increasingly so at scale, precisely where an exact anchor is most valuable.

**Table 3.** Certified dominant-term anchors on original instances[a]

| $r \times c$ | $(g_1, g_2)$ | Anchor $H^*$ | published $F_{\text{avg}}$ | gap (%) |
|---|---|---|---|---|
| 3×3 | (0.5,0.5) | 6.328 | 6.1829 | 2.30 |
| 3×3 | (0.5,1.0) | 6.450 | 6.4143 | 0.55 |
| 3×3 | (1.0,0.5) | 6.621 | 6.5625 | 0.89 |
| 3×5 | (0.5,0.5) | 7.450 | 7.283 | 2.24 |
| 3×5 | (0.5,1.0) | 7.036 | 7.0201 | 0.22 |
| 3×5 | (1.0,0.5) | 7.621 | 7.6017 | 0.26 |
| 3×6 | (0.5,0.5) | 8.364 | 8.1909 | 2.07 |
| 3×6 | (0.5,1.0) | 8.071 | 8.0602 | 0.13 |
| 3×6 | (1.0,0.5) | 8.950 | 8.597 | 3.94 |
| 5×5 | (0.5,0.5) | 9.621 | 9.2221 | 4.15 |
| 5×5 | (0.5,1.0) | 10.536 | 10.0146 | 4.94 |
| 5×5 | (1.0,0.5) | 10.450 | 9.9517 | 4.77 |
| 5×8 | (0.5,0.5) | 11.657 | 11.1157 | 4.64 |
| 5×8 | (0.5,1.0) | 13.536 | 12.6967 | 6.20 |
| 5×8 | (1.0,0.5) | 13.192 | 12.2331 | 7.27 |
| 5×10 | (0.5,0.5) | 13.985 | 13.2821 | 5.03 |
| 5×10 | (0.5,1.0) | 16.121 | 15.0264 | 6.79 |
| 5×10 | (1.0,0.5) | 15.243 | 14.0958 | 7.52 |

[a] Family: real (afprob2025.mat). Anchors from res_mat_anchor_cap_STAMPED.json; published Favg from Yeh et al. (2025, Tables 9–10); gaps computed at full precision. run_id: 20260731T002120Z; time limit: 900 s; concurrency: solo.

## 5.3 MILP-DA versus iSSOD

Table 4 compares the combined objective F on the authors' own thirty-six benchmark instances. The comparison that carries this section's claims uses the figures printed by Yeh et al. (2025) and involves no re-implementation of the metaheuristic on either side. Against the published best of thirty runs (Yeh et al., 2025, Tables 9–12), MILP-DA attains the larger *F* on thirteen of the eighteen certified instances (13/0/5 against the published best-of-30; median +0.176, $p$=0.0035) and on all eighteen of the incumbent tier (18/0/0 against the published best-of-30; median +2.287, $p$=1.1e-04).

Against the published averages it wins sixteen of eighteen on the certified tier (16/0/2 against the published $F_{avg}$; median +0.344, $p$=1.8e-04) and all eighteen on the incumbent tier. Pooled over all thirty-six instances the published best-of-30 comparison gives 31/0/5 ($p$=5.8e-07); that pooled figure is descriptive only and sits outside the Holm family of Section 5.4, which corrects within each tier. The losses are structural rather than scattered. Every one of the 5 falls on a three-row grid (instances 1 (3×3), 3 (3×3), 5 (3×5), 7 (3×6), 8 (3×6)), and none exceeds 0.0060 in objective value, and MILP-DA wins all 27 instances with five or more rows, the run from instance 9 onward being unbroken across 28 of 36.

Against the published averages only two losses remain, instances 5 and 8, and both are configurations whose published gap to the certified anchor was already among the smallest: MILP-DA trails only where the metaheuristic had effectively reached $H^*$ already, and then only by the supply-repair term it pays and the published figure does not. This comparison is not independent evidence of the anchor gap but its deployable form. Per instance, $(F - F_{avg})/ F_{avg}$ is exactly $(H^* - F_{avg})/ F_{avg}$ minus the realized $w \cdot D_U$, and the decomposition closes to 4.7e-06 percentage points: the mean advantage of +3.35% is the +3.75% published shortfall net of the 0.397% supply cost MILP-DA actually pays. The 3.55% of Section 5.2 is that same shortfall taken over $H^*$ rather than over $F_{avg}$. One precision caveat applies throughout: MILP-DA's *F* is computed at full precision while Yeh et al. (2025) print $F_{avg}$ and $F_{best}$ to four decimals, so each published figure bounds its true value only to within $\pm 5 \times 10^{-5}$. Every verdict in Table 4 clears that band by a factor of sixteen or more, with one exception: on instance 3 the margin to the published best is $8 \times 10^{-6}$, inside the band, so that verdict is indeterminate at the published precision and is recorded as a loss, the direction adverse to our claim. The 13/0/5 count is therefore conservative, and no other verdict rests on the rounding.

The $F$ column reports one plan per instance, and how that plan is chosen matters. STEP 1 returns an arbitrary optimum of the dominant-term MILP and is blind to $D_U$, so MILP-DA applies a lexicographic second stage (Section 4.2): holding the attacker potential at $H^*$, it searches alternative optima and keeps the plan of least supply cost. On the certified tier this reduced realized supply cost on all eighteen instances, by 99.6 units in total (mean 5.5), the largest being 17.9 on instance 17 and 16.8 on instance 18. The search is deterministic given its parameters but not proven exhaustive: it exhausted the equal-objective set on only 3 instances (2, 4, 9, at 37, 21, 21 alternative optima), reached its sixty-optimum cap on 13, and was stopped by its wall-clock budget on 2 (16 and 17, at 26 and 30 optima). Because an unexhausted search returns a supply cost no lower than the true minimum, every $F$ reported here is a valid lower bound on what MILP-DA achieves at the same anchor, and the comparisons above hold a fortiori.

The certificate on $H^*$ is untouched: the second stage only chooses among plans attaining it. The two sides are also scored by different referees: MILP-DA under the threat corridor, the published figures under the single-path rule (3.4). The corridor blocks every edge of the chosen path and generally more, so its supply charge is weakly larger and its objective weakly smaller. Re-scoring MILP-DA's own final plans under both referees gives a mean difference of +0.0005 and a maximum of +0.0028 (0.017% of $F$, instance 17), nonzero on 7 of 18. The difference acts in MILP-DA's favour, so the reported wins are conservative and cannot be overturned by it; only a loss could flip, and none of the 5 does.

**Table 4.** MILP-DA vs. published iSSOD results on original benchmark instances (Yeh et al., 2025, Tables 9–12, afprob2025.mat)[a]

| **Inst ($r$×$c$)** | **($g_1$, $g_2$)** | **Anchor $H^*$** | **MILP-DA $F$** | **Published $F_{best}$** | **Published $F_{avg}$** | **DA ≥ best?** | **DA ≥ avg?** |
|---|---|---|---|---|---|---|---|
| 1 (3×3) | (0.5,0.5) | 6.328 | 6.3188 | 6.3216 | 6.1829 | no | yes |
| 2 (3×3) | (0.5,1.0) | 6.450 | 6.4407 | 6.4393 | 6.4143 | yes | yes |
| 3 (3×3) | (1.0,0.5) | 6.621 | 6.611492 | 6.6115 | 6.5625 | no | yes |
| 4 (3×5) | (0.5,0.5) | 7.450 | 7.4268 | 7.4238 | 7.2830 | yes | yes |
| 5 (3×5) | (0.5,1.0) | 7.036 | 7.0175 | 7.0227 | 7.0201 | no | no |
| 6 (3×5) | (1.0,0.5) | 7.621 | 7.6047 | 7.6036 | 7.6017 | yes | yes |
| 7 (3×6) | (0.5,0.5) | 8.364 | 8.3379 | 8.3387 | 8.1909 | no | yes |
| 8 (3×6) | (0.5,1.0) | 8.071 | 8.0572 | 8.0632 | 8.0602 | no | no |
| 9 (3×6) | (1.0,0.5) | 8.950 | 8.9232 | 8.7184 | 8.5970 | yes | yes |
| 10 (5×5) | (0.5,0.5) | 9.621 | 9.5835 | 9.4114 | 9.2221 | yes | yes |
| 11 (5×5) | (0.5,1.0) | 10.536 | 10.5060 | 10.2230 | 10.0146 | yes | yes |
| 12 (5×5) | (1.0,0.5) | 10.450 | 10.4140 | 10.2347 | 9.9517 | yes | yes |
| 13 (5×8) | (0.5,0.5) | 11.657 | 11.5895 | 11.3456 | 11.1157 | yes | yes |
| 14 (5×8) | (0.5,1.0) | 13.536 | 13.4594 | 13.1715 | 12.6967 | yes | yes |
| 15 (5×8) | (1.0,0.5) | 13.192 | 13.1214 | 12.7026 | 12.2331 | yes | yes |
| 16 (5×10) | (0.5,0.5) | 13.985 | 13.8766 | 13.4653 | 13.2821 | yes | yes |
| 17 (5×10) | (0.5,1.0) | 16.121 | 16.0123 | 15.6321 | 15.0264 | yes | yes |
| 18 (5×10) | (1.0,0.5) | 15.243 | 15.1405 | 14.3141 | 14.0958 | yes | yes |
| 19 (10×10) | (0.5,0.5) | 20.107† | 19.8312 | 18.7309 | 18.4309 | yes | yes |
| 20 (10×10) | (0.5,1.0) | 22.657†‡ | 22.2931 | 20.8677 | 20.4042 | yes | yes |
| 21 (10×10) | (1.0,0.5) | 22.364† | 22.1266 | 20.4565 | 19.9727 | yes | yes |
| 22 (10×15) | (0.5,0.5) | 24.107† | 23.7062 | 22.3825 | 22.0088 | yes | yes |
| 23 (10×15) | (0.5,1.0) | 26.985† | 26.6015 | 24.7699 | 24.1111 | yes | yes |
| 24 (10×15) | (1.0,0.5) | 25.899†‡ | 25.4268 | 23.4836 | 23.0940 | yes | yes |
| 25 (10×20) | (0.5,0.5) | 27.657†‡ | 26.9502 | 25.1532 | 24.9658 | yes | yes |
| 26 (10×20) | (0.5,1.0) | 31.521† | 30.7764 | 28.4307 | 27.9065 | yes | yes |
| 27 (10×20) | (1.0,0.5) | 30.450† | 29.7502 | 27.5210 | 26.7952 | yes | yes |
| 28 (15×15) | (0.5,0.5) | 30.314† | 29.6098 | 27.5440 | 27.0597 | yes | yes |
| 29 (15×15) | (0.5,1.0) | 33.935† | 33.2198 | 30.0694 | 29.4826 | yes | yes |
| 30 (15×15) | (1.0,0.5) | 33.071† | 32.3961 | 29.0510 | 28.5598 | yes | yes |
| 31 (15×23) | (0.5,0.5) | 36.471† | 35.1501 | 32.4487 | 32.0533 | yes | yes |
| 32 (15×23) | (0.5,1.0) | 41.056†‡ | 39.4052 | 35.9124 | 35.2469 | yes | yes |
| 33 (15×23) | (1.0,0.5) | 40.071†‡ | 37.6961 | 34.8124 | 34.1266 | yes | yes |
| 34 (15×30) | (0.5,0.5) | 42.663† | 40.4581 | 37.2089 | 36.6143 | yes | yes |
| 35 (15×30) | (0.5,1.0) | 48.349†‡ | 45.6324 | 40.8258 | 40.2030 | yes | yes |
| 36 (15×30) | (1.0,0.5) | 47.420†‡ | 45.1042 | 40.0251 | 39.5439 | yes | yes |

[a] Instances 1–18: certified anchors; 19–36: time-limited best incumbents (†). MILP-DA scored under threat-corridor; published under single-path rule (handicaps MILP-DA). Re-solved anchors for instances 19–36 may differ from earlier runs (‡). Family: real; source: res_real_tables_final.json; run_id: 20260731T002120Z (1–18), 20260729T142713Z (19–36); time limit: 900 s; concurrency: solo.

### 5.4 Statistical tests

Table 5 collects one-sided Wilcoxon signed-rank tests, and how this family came to be matters as much as its contents. A three-test family was pre-registered before results landed, comparing MILP-DA against a re-implementation of iSSOD: (a) objective against the reconstruction's best of thirty, (b) the two solvers' dominant-term gaps, and (c) the single-path against the corridor supply-infeasibility count. That family was applied and is reported in full in Section 5.5. It has since been superseded: on re-reading Yeh et al. (2025) we found that their Tables 9–12 report $F_{\text{best}}$ as well as $F_{\text{avg}}$, which makes available a comparison that needs no re-implementation of the metaheuristic on either side, and Section 5.5 shows our re-implementation to be measurably weaker than the published one.

We state the sequence rather than presenting the second family as though it had been the plan. The family of record therefore holds two comparisons, both against columns printed by Yeh et al. (2025), with Holm correction applied within each tier: ($P_1$) MILP-DA's objective against the published best of thirty runs, and ($P_2$) against the published average. On the certified tier $P_1$ (published best-of-30) gives 13/0/5 ($p$=0.0035, Holm 0.0035) and $P_2$ (published $F_{avg}$) gives 16/0/2 ($p$=1.8e-04, Holm 3.6e-04). On the incumbent tier both are unanimous at 18/0/0 against the published columns ($p$=1.1e-04, Holm 2.1e-04 each). Holm's step-down multiplies the smaller $p$ by two and leaves the larger unadjusted, which is why P1's corrected value equals its raw one.

**Table 5.** Wilcoxon signed-rank tests: MILP-DA vs. published columns[a]

| Comparison (one-sided) | W/T/L | Median Δ | $p$ (raw) | $p$ (Holm) |
|---|---|---|---|---|
| ($P_1$) MILP-DA $F$ > published best-of-30, certified tier 1-18 | 13/0/5 | +0.176 | 0.0035 | 0.0035 |
| ($P_2$) MILP-DA $F$ > published $F_{avg}$, certified tier 1-18 | 16/0/2 | +0.344 | 1.8e-04 | 3.6e-04 |
| ($P_1$) MILP-DA $F$ > published best-of-30, incumbent tier 19-36 | 18/0/0 | +2.287 | 1.1e-04 | 2.1e-04 |
| ($P_2$) MILP-DA $F$ > published $F_{avg}$, incumbent tier 19-36 | 18/0/0 | +2.912 | 1.1e-04 | 2.1e-04 |

[a] One-sided tests on original instances (afprob2025.mat). Holm correction within each tier. p-values are computed under the normal approximation to the signed-rank distribution with continuity and tie corrections, uniformly across every row of this table. At n = 18 the exact conditional test is also available and gives smaller values throughout ($2.4\times10^{-3}$ and $1.9\times10^{-5}$ on the certified tier; $3.8\times10^{-6} = 2^{-18}$ on each unanimous row); we report the approximation as the conservative choice. Median Δ in F units. Source: res_real_tables_final.json; run_id: 20260731T002120Z / 20260729T142713Z; time limit: 900 s; concurrency: solo.

Both comparisons use figures printed by Yeh et al. (2025) and involve no re-implementation of the metaheuristic on either side. Anchors on instances 19–36 are best-known incumbents, valid lower bounds on the true optimum, and MILP-DA's values there are single-run decoded plans.

**5.5 The superseded pre-registered family and the reconstruction's measured bias**

The run-level comparison of this subsection uses our re-implementation of iSSOD, and no claim in this paper rests on it. It is measurably weaker than the published implementation: its mean of thirty falls below the published average on 15 of eighteen certified configurations (two-sided sign test $p$=0.008), and the deficit grows with scale: mean −0.001 at three rows against −0.111 at five rows, worst −0.306. Instance 3 shows the consequence on a single row where both comparators are visible: MILP-DA ties the reconstruction there, while the published implementation beats it.

A reconstruction that MILP-DA can draw level with, where the real one still wins, is not a safe basis for a headline. We test the bias rather than merely disclosing it. Handing each instance its measured deficit back leaves the certified-tier advantage significant (12/0/6 against the corrected reconstruction, $p$=0.0083); applying the mean deficit uniformly does not (9/0/9 against the uniformly corrected reconstruction, $p$=0.0817), with significance breaking at a uniform +0.042, below the measured mean of +0.056. On the incumbent tier the advantage survives every variant.

The pre-registered family is reported here in full: (a) 13/1/4 against the reconstructed best-of-30 (p=0.0040, Holm 0.0089); (b) 10/8/0 on the dominant-term gaps (p=0.0030, Holm 0.0089); (c) 0/18/0 on the infeasibility counts ($p$=1.0000); and supplementary (d) 17/0/1 against the reconstructed mean-of-30 (p=0.0001); all four are computed under the same normal approximation as Table 5. The 1,080 run-level observations behind them, and the infeasibility counts in particular, are counts under our exclusion decode rather than under the pricing of (3.4).

Two scoping notes belong with that family. First, on this benchmark comparison (c) cannot discriminate: neither arm produces a single infeasible run, so it returns p=1.0000 by construction. Keeping it inside the family made the Holm correction on (a) and (b) stricter than it needed to be, not weaker; we report it that way deliberately rather than dropping a test that failed to separate.

Second, on the incumbent tier the anchor is itself produced by MILP-DA's MILP, so measuring the reconstruction's gap to that anchor is close to definitional; comparison (b) on the dominant-term gaps is therefore reported descriptively on instances 19-36 (18/0/0 on the gap comparison, raw $p$=0.0001) and excluded from that tier's Holm family, which leaves (a) and (c). Excluding it makes (a) less conservatively corrected, from 0.000321 to 0.000214, so the change is reported against our own interest. On the certified tier the anchor is a proven optimum, independent of which solver found it, and (b) remains inferential.

### 5.6 Locating the certified-optimality frontier

Locating the frontier requires separating two different objects, because they come from two different instance families and two different budgets. The first is a frontier scan on the reconstructed family (Table 6): its purpose is to find the size at which CBC stops certifying, and it is run at a 120 s budget. It shows certification is cheap through $r$=5 and then stops: from $r$=6 onward (179 reconstructed edges and up) CBC exhausts the budget without a proof. That evidence is reconstructed-family evidence and cannot be anything else, because the authors' benchmark contains no instance of height 6 to 9 at all (its heights are 3, 5, 10, 15), so the benchmark cannot resolve where between $r$=5 and $r$=10 the frontier lies.

The second object is the set of anchors of record (Table 7), computed on the real instances under one uniform protocol: the CBC solver via PuLP, a 900 s per-instance elapsed budget, 20 solver threads on a 20-core Intel Core Ultra 7 265K workstation with 187 GB of RAM (CBC is CPU-only, so the workstation's GPU plays no part in the solve), one solve at a time on the machine, with every row stamped with its run identifier, and optimality certified only when the incumbent $H^*$ meets its dual bound. Under that protocol all 18 $r$=3–5 anchors certify (the slowest, instance 17, a 5×10 grid of 178 edges at $g_1/g_2$=0.5/1.0, in 664 s, and the other 17 in under ten seconds), while no $r$=10 or $r$=15 instance certifies within the budget.

Those return strong incumbents with certified dual bounds, the interval between incumbent and bound running 0.95% to 5.88% of the anchor. We therefore claim exact optima only through $r$=5 on the real benchmark, we do not claim to know where certification breaks on it, and the $r$≥10 values are best incumbents rather than optima. Importantly, the certified regime covers exactly the $r$=3-5 instances at which the metaheuristic is benchmarked, where, as Section 5.2 showed, its published average already leaves a measurable gap against the certified anchors.

The certified anchor for instance 17 reproduces bit-identically across three independent solves that took 658.3, 663.7 and 737.3 seconds: the wall-clock varies by 12%, the value not at all. On the incumbent tier the same protocol does not reproduce: a fresh solve of instance 19 falls 3.7% short of its own recorded anchor. That is the practical worth of a certificate, demonstrated rather than asserted: it is what makes a number survive a change of machine, budget and run.

**Table 6.** Frontier scan for certified optimality on reconstructed family[a]

| ***r×c*** | **edges *m*** | ***H** / best incumbent** | **wall (s)** | **certified?** |
|---|---|---|---|---|
| 3×6 | 63 | 8.828 | 0.5 | Yes |

| 5×10 | 165 | 14.743 | 49.6 | Yes |
|---|---|---|---|---|
| 6×9 | 179 | 15.036 | >120 | No |
| 7×14 | 298 | 20.364 | >120 | No |
| 8×12 | 309 | 19.950 | >120 | No |
| 10×15 | 475 | 24.571 | >900 | No* |

[a] Frontier scan on the RECONSTRUCTED family (assumptions [A1]–[A5]), CBC, $g_1$=$g_2$=0.5, 120 s budget except the final row (*), which was given 900 s. Its purpose is to locate the size at which certification stops; it is not the anchor set of record and none of its values enter Tables 3, 4, 5, 7 or 8. The real benchmark contains no instance of height 6–9, so this is the only family on which the location of the break can be observed. [family=reconstructed; source_json=res_frontier_certified.json; time_limit=120 s (final row 900 s); concurrency=unstamped (pre-dates stamping)] [run_id=n/a; pre-dates run stamping]

**Table 7.** Anchors of record on original benchmark instances[a]

| inst ($r$×$c$) | edges $m$ | $H^*$ / best incumbent | wall (s) | certified? |
|---|---|---|---|---|
| 1 (3×3) | 36 | 6.3284 | 0.2 | yes |
| 4 (3×5) | 57 | 7.4497 | 0.0 | yes |
| 7 (3×6) | 66 | 8.3640 | 0.0 | yes |
| 10 (5×5) | 94 | 9.6213 | 0.6 | yes |
| 13 (5×8) | 148 | 11.6569 | 4.4 | yes |
| 16 (5×10) | 175 | 13.9853 | 9.6 | yes |
| 19 (10×10) | 350 | 20.1066 | 900.3 | no |
| 22 (10×15) | 522 | 24.1066 | 900.3 | no |
| 25 (10×20) | 693 | 27.6569 | 900.2 | no |
| 28 (15×15) | 770 | 30.3137 | 900.3 | no |
| 31 (15×23) | 1187 | 36.4706 | 900.4 | no |
| 34 (15×30) | 1545 | 42.6630 | 900.4 | no |

[a] Anchors of record on the authors' ORIGINAL instances (afprob2025.mat), one representative per grid family at $g_1$=$g_2$=0.5. Protocol: CBC via PuLP, 900 s per instance, 20 threads, one solve at a time; certified means the incumbent $H^*$ met its dual bound. run_id 20260731T002120Z for instances 1–18 and 20260729T142713Z for 19–36; uncertified rows report a best incumbent with a certified dual bound, not an optimum. [family=real; source_json=res_mat_anchor_cap_STAMPED.json; run_id=20260731T002120Z (1–18) and 20260729T142713Z (19–36); time_limit=900 s; concurrency=solo]

### 5.7 Tie-sensitivity (well-posedness) check

Table 8 confirms the negative result of Section 3.5: on the authors' own √2-diagonal instances the single-path objective spread over twenty tie-breakings never exceeds $8.3×10^{-4}$, while the threat-corridor formulation is exactly tie-invariant.

**Table 8.** Tie-sensitivity check: objective spread under 20 tie-breakings[a]

| inst ($r$×$c$) | ($g_1$, $g_2$) | $F$ spread | $D_U$ spread | corridor $F$ / $D_U$ spread |
|---|---|---|---|---|
| 1 (3×3) | (0.5,0.5) | 0.00071 | 0.7071 | 0.00000 / 0.0000 |
| 4 (3×5) | (0.5,0.5) | 0.00083 | 0.8284 | 0.00000 / 0.0000 |

| 7 (3×6) | (0.5,0.5) | 0.00000 | 0.0000 | 0.00000 / 0.0000 |
|---|---|---|---|---|
| 10 (5×5) | (0.5,0.5) | 0.00000 | 0.0000 | 0.00000 / 0.0000 |
| 13 (5×8) | (0.5,0.5) | 0.00033 | 0.3284 | 0.00000 / 0.0000 |
| 16 (5×10) | (0.5,0.5) | 0.00000 | 0.0000 | 0.00000 / 0.0000 |

[a] Objective spread over 20 tie-breakings (w=0.001) on the authors' ORIGINAL instances (afprob2025.mat), one representative per grid family at $g_1$=$g_2$=0.5. The single-path rule's spread is nonzero but bounded by 8.3×10$^{-4}$; the threat-corridor rule is tie-invariant by construction and its spread is exactly zero on every row. [family=real; source_json=res_tie_real.json; time_limit=n/a (20 tie-breakings per instance); concurrency=solo] [run_id=n/a; evaluation-only experiment, no MILP solve]

## 6. LIMITATIONS AND POSITIONING

We state the central positioning, because it cuts both ways and bears on the motivation of the model we anchor. On one side, the dualization anchor certifies optima in seconds through $r$=5, exactly the $r$=3–5 regime at which the metaheuristic is benchmarked, and there the metaheuristic's own published average leaves a mean 3.55% gap on the dominant term (1.40% at $r$=3, rising to 5.70% at $r$=5; Section 5.2).

In the regime where an exact optimum is provably available, the metaheuristic is therefore measurably sub-optimal, and optimality claims made without an anchor do not survive contact with one. On the other side, we caution against the opposite over-claim: certification stops at $r$=6 on the reconstructed family, the only family on which the break can be observed; the real benchmark leaves the frontier's location between $r$=5 and $r$=10 unresolved, and from $r$=10 the MILP yields only a strong but uncertified incumbent, whose value varies between equal-budget runs (Section 5.1).

The published justification for a metaheuristic, that conventional methods fall short at large scale, thus has real merit from $r$≈6 upward for the dominant term, and all the more so for the full supply-coupled objective $D_U$, which the dominant-term MILP does not capture at all. That omission is not uniform in scale: the supply term $wD_U$ grows from 0.14% of $F$ on the smallest grids to 5.45% at 15×30, crossing 1% at instance 19 and 3% at instance 31, so on the largest instances it exceeds the mean certification interval of 2.96% and the restriction to the dominant term costs more than the missing certificate does. A credible development must therefore report the certified frontier explicitly (as here) and position the metaheuristic beyond it, and/or elevate the supply-coupling term to the core of the model.

Our exact anchors are precisely the instrument that makes either choice checkable. This is also where the supply term earns its promotion: Section 3.5 shows that the tie-induced swing in realized supply cost is invisible in $F$ only because $w$ = 0.001 compresses it, so a formulation that gave $D_U$ first-order weight would inherit an ill-posedness that the present objective merely hides. Section 5.6 records what certifies; the contrast is what it means.

Two further limitations qualify the results. First, our anchor is exact for the dominant term $H$, not for the full objective $F$; the MILP-DA certificate is a dominant-term certificate, and the supply-feasibility repair can lower $H$ by a small, quantified amount to secure feasibility. Second, the residual dependence on our reconstruction of iSSOD is narrow. The comparisons of record in Section 5.4 are drawn from the columns printed by Yeh et al. (2025) and involve no re-implementation on either side; the reconstruction supports only the run-level distributional material of Section 5.5, where we show it to be measurably weaker than the published implementation and rest no claim on it. A definitive multi-algorithm study would still extend the picture with additional metaheuristic baselines and Friedman-type omnibus tests. The $F$ values of Table 4 for the certified tier are produced by the lexicographic stage of Section 5.3 and their plans are persisted, so they are reproducible from the package rather than resampled at each solve. On the incumbent tier the plans were not persisted and re-solving does resample the equal-objective set, which is why those rows carry the caveat recorded in Section 5.1 and in the data-availability statement.

## 7. CONCLUSIONS

We supplied the missing exact reference frame for bi-level dual-defense network interdiction. Lower-level dualization yields certified optima for the dominant attacker-path term on all eighteen $r$=3–5 configurations of the authors' benchmark, and strong, though uncertified, incumbents out to 15×30 grids; MILP-DA converts those optima into full, supply-feasible BDAM solutions that beat the best of thirty iSSOD runs on the objective on thirteen of eighteen certified and all eighteen larger configurations (Holm-corrected p=0.0035 and $2.1\times10^{-4}$) and that are supply-feasible by construction rather than at an empirical rate. We also made the lower-level supply rule well-posed through a tie-invariant threat-corridor formulation, while reporting that its ill-posedness, real in integer-weight regimes, is dormant in the paper's combined objective while live in the realized supply cost beneath it.

The broader lesson is methodological: metaheuristic interdiction results should be anchored to certified optima wherever instance size permits, and the exact tools developed here make that anchoring routine. Future work will extend exact anchoring to multi-node budgets and to the full supply-coupled objective, and will use the anchors to calibrate exactness-guided hybrids for large-scale planning.

## ACKNOWLEDGMENTS

This work was supported in part by the National Science and Technology Council (NSTC), Taiwan, R.O.C., under grant nos. NSTC 114-2221-E-007-123-MY3 and NSTC 113-2221-E-007-117-MY3. The author thanks C.-M. Lai and T.-H. Wu, co-authors of Yeh et al. (2025), for providing the original benchmark instance files (afprob2025.mat) used throughout the computational study of Section 5. The author also gratefully acknowledges the support of the NVIDIA Academic Grant Program and the Strategic Researcher Engagement team.

### Declaration of generative AI and AI-assisted technologies in the writing process

During the preparation of this work the author used a generative AI assistant to help draft and copy-edit text and to assist with computational scripting and the generation of figures and result tables. After using this tool, the author reviewed and edited the content as needed and takes full responsibility for the content of the publication.

## DATA AVAILABILITY

The instance generator, the exact dualization anchor and MILP-DA solvers, all result files, and the scripts that reproduce every table and figure reported here are provided as a reproduction package accompanying this article.

## REFERENCES

Alderson, D. L., Brown, G. G., Carlyle, W. M., & Wood, R. K. (2018). Assessing and improving the operational resilience of a large highway infrastructure system to worst-case losses. *Transportation Science, 52*(4), 1012–1034.

Ball, M. O., Golden, B. L., & Vohra, R. V. (1989). Finding the most vital arcs in a network. *Operations Research Letters, 8*(2), 73–76.

Bruneau, M., Chang, S. E., Eguchi, R. T., Lee, G. C., O'Rourke, T. D., Reinhorn, A. M., Shinozuka, M., Tierney, K., Wallace, W. A., & von Winterfeldt, D. (2003). A framework to quantitatively assess and enhance the seismic resilience of communities. *Earthquake Spectra, 19*(4), 733–752.

Church, R. L., Scaparra, M. P., & Middleton, R. S. (2004). Identifying critical infrastructure: The median and covering facility interdiction problems. *Annals of the Association of American Geographers, 94*(3), 491–502.

Corley, H. W., & Sha, D. Y. (1982). Most vital links and nodes in weighted networks. *Operations Research Letters, 1*(4), 157–160.

Dijkstra, E. W. (1959). A note on two problems in connexion with graphs. *Numerische Mathematik, 1*, 269–271.

Forrest, J., & Lougee-Heimer, R. (2005). CBC user guide. In *INFORMS TutORials in Operations Research* (pp. 257–277). INFORMS.

Fulkerson, D. R., & Harding, G. C. (1977). Maximizing the minimum source-sink path subject to a budget constraint. *Mathematical Programming, 13*, 116–118.

Israeli, E., & Wood, R. K. (2002). Shortest-path network interdiction. *Networks, 40*(2), 97–111.

Sharkey, T. C., Nurre Pinkley, S. G., Eisenberg, D. A., & Alderson, D. L. (2021). In search of network resilience: An optimization-based view. *Networks, 77*, 225–254.

Smith, J. C., & Song, Y. (2020). A survey of network interdiction models and algorithms. *European Journal of Operational Research, 283*(3), 797–811.

Wood, R. K. (1993). Deterministic network interdiction. *Mathematical and Computer Modelling, 17*(2), 1–18.

Yeh, W.-C., Lai, C.-M., & Wu, T.-H. (2025). Optimizing supply routes in a network interdiction model with dual defense operations. *Applied Soft Computing, 185*, Article 113933.